# Buildings, Group Extensions and the Cohomology of Congruence Subgroups

**Alejandro Adem***
Mathematics Department
University of Wisconsin
Madison WI 53706

## Introduction

Let $T_n$ denote the Tits Building associated to the discrete group $\mathrm{SL}_n(\mathbb{Z})$. In this note we will be interested in computing the equivariant cohomology $H^*_\Gamma(|T_n|, \mathbb{Z})$, where $|T_n|$ is the geometric realization of $T_n$, and $\Gamma$ is a torsion-free discrete subgroup of finite index in $\mathrm{SL}_n(\mathbb{Z})$. In the special case $n = 3$, we will see that this equivariant cohomology in fact corresponds to the ordinary cohomology of a 4-dimensional Poincaré Duality group which has $\Gamma$ as its quotient.

This is then applied to partially compute the cohomology of the level $p$ congruence subgroups ($p$ an odd prime) in $\mathrm{SL}_3(\mathbb{Z})$. In particular we obtain

**Theorem**
*Let $\Gamma(p)$ denote the level $p$ congruence subgroup of $SL_3(\mathbb{Z})$; then*

$$dim_{\mathbb{Q}}\ H^3(\Gamma(p), \mathbb{Q}) \geq \frac{1}{12}(p^3 - 1)(p^3 - 3p^2 - p + 15) + 1.$$

It turns out that this method can be used for other rank two groups such as $\mathrm{Sp}_4(\mathbb{Z})$ and $\mathrm{G}_2(\mathbb{Z})$. In §6 we apply this to the level $p$ congruence subgroups in $\mathrm{Sp}_4(\mathbb{Z})$, obtaining the cohomology of the relevant parabolic subgroups and from this a lower bound for the fourth betti number:

**Theorem**
*Let $\Gamma(p)$ denote the level $p$ congruence subgroup in $Sp_4(\mathbb{Z})$; then*

$$dim_{\mathbb{Q}}\ H^4(\Gamma(p), \mathbb{Q}) \geq \frac{1}{24}(p^4 - 1)(2p^3 - 3p^2 - 2p + 27) + 1.$$

Our methods are mostly algebraic, involving techniques from group extensions and their cohomology. We offer a systematic approach which works with general coefficients and in addition we provide complete information on the cohomology of the parabolic subgroups.

* Partially supported by an NSF Grant

One of the key facts which we prove is that the congruence subgroups are homomorphic images of Poincaré duality groups with kernel a free, infinitely generated group. These groups are interesting in their own right, and probably deserve further attention. They can be realized as the fundamental groups of the Borel–Serre compactification; this is discussed in §7.

The application to congruence subgroups in $\mathrm{SL}_3(\mathbb{Z})$ is contained in the paper by Lee and Schwermer (see [LS]), although it is perhaps not as well known as it should be. In fact most of the work in this paper was done without knowledge of the existence of this formidable piece of work! We offer a somewhat simplified and purely cohomological approach, providing an elementary and explicit computation. On the other hand, the application to congruence subgroups in $\mathrm{Sp}_4(\mathbb{Z})$ seems to be new.

In some sense this is an expository paper. One of its goals is to put computations for congruence subgroups within the framework of group cohomology. Complete answers involve assembling a wide array of mathematical techniques beyond our scope, but at the very least one hopes this paper will motivate the reader to learn the more number–theoretic aspects of the subject (see [A], [AGG]).

The author is grateful to A. Ash and L. Solomon for helpful comments.

## §0. Preliminaries

We recall that $\mathrm{SL}_n(\mathbb{Z})$ is a group of virtual cohomological dimension equal to $N = \frac{n(n-1)}{2}$. The module known as the Steinberg module St is a dualizing module for any torsion-free subgroup $\Gamma$ of finite index in $\mathrm{SL}_n(\mathbb{Z})$; that is there are isomorphisms

$$H^i(\Gamma, M) \cong H_{N-i}(\Gamma, \mathrm{St} \otimes M)$$

where $M$ is any $\mathbb{Z}\Gamma$ module, $i$ any integer and these equivalences are given by taking cap product with a fixed element $z \in H_N(\Gamma, \mathrm{St})$; dually, we have isomorphisms

$$H_i(\Gamma, M) \cong H^{N-i}(\Gamma, \mathrm{Hom}(\mathrm{St}, M)).$$

Recall that $\mathrm{SL}_n(\mathbb{Q})$ acts naturally on the set of all proper non-trivial subspaces of $\mathbb{Q}^n$. This defines a partially ordered set $T_n$ and hence its geometric realization $|T_n|$ which can be thought of as an $\mathrm{SL}_n(\mathbb{Z})$-CW complex of dimension $n-2$. Furthermore we have

$$H_*(|T_n|, \mathbb{Z}) \cong \begin{cases} \mathbb{Z} & * = 0 \\ \mathrm{St} & * = n-2 \\ 0 & \text{otherwise.} \end{cases}$$

For details on these facts we refer to [Br1] and [Br2].

Next we recall the notion of equivariant cohomology: if $G$ is a group acting on a space $X$, then

$$H^*_G(X) = H^*(X \times_G EG, \mathbb{Z})$$

where $EG$ is the universal free contractible $G$ space and $X \times_G EG = X \times EG/G$ (diagonal action).

We also recall that if $X$ is a G-CW complex, then there exist two spectral sequences for computing $H^*_G(X)$, namely
(I) $E_2^{p,q} = H^p(G, H^q(X, \mathbb{Z})) \Rightarrow H_G^{p+q}(X)$
and
(II) $E_1^{p,q} = H^q(G, C^p(X)) \Rightarrow H_G^{p+q}(X)$
where $C^*(X)$ denotes the module of cellular $p$-cochains on $X$.

We remark here for later use that if

$$C_p(X) \cong \bigoplus_{i \in I} \mathbb{Z}[G/G_{\sigma_i}],$$

then

$$H^q(G, C^p(X)) \cong \bigoplus_{i \in I} H^q(G_{\sigma_i}, \mathbb{Z}).$$

## §1. The Equivariant Cohomology of $|T_n|$.

We will prove the following

**Theorem 1.1**. *Let $\Gamma \subseteq SL_n(\mathbb{Z})$ denote a torsion free subgroup of finite index, then there is a duality isomorphism*

$$H^i_\Gamma(|T_n|) \cong H^\Gamma_{\frac{n^2+n-4}{2}-i}(|T_n|).$$

**Proof**. To begin we observe that there is a fibration

$$\begin{array}{ccc} |T_n| & \longrightarrow & |T_n| \times_\Gamma E\Gamma \\ & & \downarrow \\ & & B\Gamma \end{array}$$

where $\dim|T_n| = n-2$, $\dim B\Gamma = N = \frac{n(n-1)}{2}$ and hence $\dim |T_n| \times_\Gamma E\Gamma = n-2+N = \frac{n^2+n-4}{2}$. Let us denote this number by $M$. Now in this case the spectral sequence of type (I) in homology degenerates to a long exact sequence

$$\cdots \to H_{i-n+2}(\Gamma, \mathrm{St}) \to H^\Gamma_i(|T_n|) \to H_i(\Gamma, \mathbb{Z}) \to H_{i-n+1}(\Gamma, \mathrm{St}) \to \cdots.$$

Note in particular that

$$H^\Gamma_M(X) \cong H_N(\Gamma, \mathrm{St}) \cong H^0(\Gamma, \mathbb{Z}) \cong \mathbb{Z}.$$

Hence there is a "fundamental class" $\mu \in H^\Gamma_M(|T_n|)$ which can be identified with the class $z \in H_N(\Gamma, \mathrm{St})$ which induces the duality isomorphisms $H^i(\Gamma, M) \cong H_{N-i}(\Gamma, \mathrm{St} \otimes M)$ via the cap product. We claim that

$$\cap\mu : H^i_\Gamma(|T_n|) \to H^\Gamma_{M-i}(|T_n|)$$

induces a duality isomorphism.

To prove this we look at the long exact sequence in cohomology:

$$\cdot \to H^i(\Gamma,\mathbb{Z}) \to H^i_\Gamma(|T_n|) \to H^{i-n+2}(\Gamma,\mathrm{St}^*) \to H^{i+1}(\Gamma,\mathbb{Z}) \to \cdots.$$

Cap product with $\mu$ induces maps from one long exact sequence to the other:

$$\begin{array}{ccccccccc}
\cdots \to & H^i(\Gamma,\mathbb{Z}) & \to & H^i_\Gamma(|T_n|) & \to & H^{i-n+2}(\Gamma,\mathrm{St}^*) & \to & H^{i+1}(\Gamma,\mathbb{Z}) & \to \cdots \\
& \downarrow \cong & & \downarrow \cap\mu & & \downarrow \cong & & \downarrow \cong & \\
\cdots \to & H_{N-i}(\Gamma,\mathrm{St}) & \to & H^\Gamma_{M-i}(|T_n|) & \to & H_{M-i}(\Gamma,\mathbb{Z}) & \to & H_{N-i-1}(\Gamma,\mathrm{St}) & \to \cdots
\end{array}$$

Applying the five-lemma completes the proof.

□

Consider the special case when $n = 3$; then $|T_3|$ is an infinite graph, with free fundamental group $F$. In this case $|T_n| \times_\Gamma E\Gamma$ is aspherical and 4 dimensional; hence if $Q = \pi_1(|T_n| \times_\Gamma E\Gamma)$ we have an extension

$$1 \to F \to Q \to \Gamma \to 1$$

where $Q$ is a 4-dimensional Poincaré Duality group with homology of *finite type.* In particular we have a long exact sequence:

$$0 \to H^1(\Gamma,\mathbb{Z}) \to H^1(Q,\mathbb{Z}) \to H^0(\Gamma,\mathrm{St}^*) \to H^2(\Gamma,\mathbb{Z}) \to H^2(Q,\mathbb{Z})$$

$$\to H^1(\Gamma,\mathrm{St}^*) \to H^3(\Gamma,\mathbb{Z}) \to H^3(Q,\mathbb{Z}) \to H^2(\Gamma,\mathrm{St}^*)$$

$$\to H^4(\Gamma,\mathbb{Z}) \to H^4(Q,\mathbb{Z}) \to H^3(\Gamma,\mathrm{St}^*) \to 0$$

which simplifies to yield

$$\begin{array}{llllll}
0 \to H^1(\Gamma,\mathbb{Z}) & \to H^1(Q,\mathbb{Z}) & \to H_3(\Gamma,\mathbb{Z}) & \to H^2(\Gamma,Z) & \to H^2(Q,\mathbb{Z}) & \to \\
& \to H_2(\Gamma,\mathbb{Z}) & \to H^3(\Gamma,\mathbb{Z}) & \to H^3(Q,\mathbb{Z}) & \to H_1(\Gamma,\mathbb{Z}) & \to 0
\end{array} \tag{1.2}$$

and $H^4(Q,\mathbb{Z}) \cong \mathbb{Z}$.

## §2. Cohomology of Parabolic Subgroups

From now on we specialize to the case $n = 3$ and $\Gamma = \Gamma(p)$, a level $p$ congruence subgroup, $p$ an odd prime and $Q = Q(p)$. We will compute the integral cohomology of its intersections with the parabolic subgroups in $\mathrm{SL}_3(\mathbb{Z})$.

To begin let $B \subseteq \mathrm{SL}_n(\mathbb{Z})$ denote the subgroup of upper triangular matrices in $\mathrm{SL}_3(\mathbb{Z})$. Then

$$B \cap \Gamma(p) = \left\{ \begin{pmatrix} 1 & c_1 & c_2 \\ 0 & 1 & c_3 \\ 0 & 0 & 1 \end{pmatrix} \ p \text{ divides } c_1,\ c_2, \text{ and } c_3 \right\}.$$

The following elements generate this group

$$a = \begin{pmatrix} 1 & p & 0 \\ 0 & 1 & 0 \\ 0 & 0 & 1 \end{pmatrix}, \qquad b = \begin{pmatrix} 1 & 0 & 0 \\ 0 & 1 & p \\ 0 & 0 & 1 \end{pmatrix}, \qquad c = \begin{pmatrix} 1 & 0 & p \\ 0 & 1 & 0 \\ 0 & 0 & 1 \end{pmatrix}$$

and $aba^{-1}b^{-1} = c^p$.

In fact we may express $B \cap \Gamma(p)$ as a central extension

$$1 \to \mathbb{Z} \to B \cap \Gamma(p) \to \mathbb{Z} \oplus \mathbb{Z} \to 1$$

where $\langle c \rangle$ is the central subgroup and $\overline{a}$, $\overline{b}$ generate the quotient. We have an exact sequence in cohomology (over $\mathbb{Z}$):

$$0 \to \mathbb{Z} \oplus \mathbb{Z} \to H^1(B \cap \Gamma(p)) \to \mathbb{Z} \xrightarrow{d_2} \mathbb{Z} \to H^2(B \cap \Gamma(p)) \to \mathbb{Z} \oplus \mathbb{Z} \to 0$$

and $H^3(B \cap \Gamma(p)) \cong \mathbb{Z}$. From the extension data, coker $d_2 \cong \mathbb{Z}/p$ and hence we have proved

**Proposition 2.1.**

$$H^*(B \cap \Gamma(p), \mathbb{Z}) \cong \begin{cases} \mathbb{Z} & * = 0, 3 \\ \mathbb{Z} \oplus \mathbb{Z} & * = 1 \\ \mathbb{Z} \oplus \mathbb{Z} \oplus \mathbb{Z}/p & * = 2 \\ 0 & \text{otherwise.} \end{cases}$$

□

This is a 3-dimensional Poincaré Duality group. Next we consider the parabolic subgroup $P_1 \subseteq \mathrm{SL}_3(\mathbb{Z})$ the subgroup with zeros in the second and third entries of its first column. We have

$$P_1 \cap \Gamma(p) = \left\{ \begin{pmatrix} 1 & x & y \\ 0 & & \\ & & A \\ 0 & & \end{pmatrix} \middle| \begin{array}{c} p \text{ divides } x \text{ and } y \\ A \in \Gamma_2(p) \end{array} \right\}.$$

Let

$$E = \left\{ \begin{pmatrix} 1 & x & y \\ 0 & 1 & 0 \\ 0 & 0 & 1 \end{pmatrix} \; p \text{ divides } x \text{ and } y \right\},$$

then $E \cong \mathbb{Z} \oplus \mathbb{Z}$ and if we embed $\Gamma_2(p)$ in $\Gamma(p)$ via $A \mapsto \begin{pmatrix} 1 & \cdots \\ \vdots & \boxed{A} \end{pmatrix}$ we obtain a semidirect product decomposition

$$\boxed{P_1 \cap \Gamma(p) = E \times_T \Gamma_2(p)}.$$

If $A = \begin{pmatrix} a & b \\ c & d \end{pmatrix} \in \Gamma_2(p)$, the action on $E$ is given by

$$\begin{pmatrix} 1 & x & y \\ 0 & 1 & 0 \\ 0 & 0 & 1 \end{pmatrix} \overset{T_A}{\mapsto} \begin{pmatrix} 1 & ax+cy & bx+dy \\ 0 & 1 & 0 \\ 0 & 0 & 1 \end{pmatrix}$$

or

$$x \mapsto ax + cy$$
$$y \mapsto bx + dy$$

and hence

$$T_A \begin{pmatrix} x \\ y \end{pmatrix} = \begin{pmatrix} a & c \\ b & d \end{pmatrix} \begin{pmatrix} x \\ y \end{pmatrix} = A^T \begin{pmatrix} x \\ y \end{pmatrix}.$$

We can now compute $H^*(P \cap \Gamma(p), \mathbb{Z})$; we use the Lyndon-Hochschild-Serre spectral sequence for the extension above. It has $E_2$ term

$$E_2^{p,q} = H^p(\Gamma_2(p), H^q(E, \mathbb{Z})).$$

Note that as $\Gamma_2(p)$ is free (and hence 1-dimensional), there are no differentials and it collapses at $E_2$. We obtain

$$H^i(P_1 \cap \Gamma(p), \mathbb{Z}) \cong \begin{cases} \mathbb{Z} & i = 0 \\ H^1(\Gamma_2(p), \mathbb{Z}) \oplus H^1(E, \mathbb{Z})^{\Gamma_2(p)} & i = 1 \\ H^1(\Gamma_2(p), H^1(E, \mathbb{Z})) \oplus H^2(E, \mathbb{Z})^{\Gamma_2(p)} & i = 2 \\ H^1(\Gamma_2(p), H^2(E), \mathbb{Z}) & i = 3. \end{cases}$$

Now the $\Gamma_2(p)$ action on $H^1(E)$ is the natural action on $M = \mathbb{Z} \oplus \mathbb{Z}$, and hence that on $H^2(E)$ is via the determinant, from which $H^1(E, \mathbb{Z})^{\Gamma_2(p)} = 0$, $H^2(E, \mathbb{Z})^{\Gamma_2(p)} \cong \mathbb{Z}$ and we have proved

**Proposition 2.2**.

$$H^i(P_1 \cap \Gamma(p), \mathbb{Z}) \cong \begin{cases} \mathbb{Z} & i = 0 \\ H^1(\Gamma_2(p), \mathbb{Z}) & i = 1 \\ H^1(\Gamma_2(p), M) \oplus \mathbb{Z} & i = 2 \\ H^1(\Gamma_2(p), \mathbb{Z}) & i = 3 \\ 0 & \text{otherwise.} \end{cases}$$

□

Denote by $M^*$ the dual 2-dimensional representation of $\Gamma_2(p)$. Then it is direct to verify that for the other parabolic subgroup $P_2$, we have

**Proposition 2.3**.

$$H^i(P_2 \cap \Gamma(p), \mathbb{Z}) \cong \begin{cases} \mathbb{Z} & i = 0 \\ H^1(\Gamma_2(p), \mathbb{Z}) & i = 1 \\ H^1(\Gamma_2(p), M^*) \oplus \mathbb{Z} & i = 2 \\ H^1(\Gamma_2(p), \mathbb{Z}) & i = 3 \\ 0 & \text{otherwise.} \end{cases}$$

□

Note that $\Gamma_2(p)$ is a free group of rank $1 + \frac{(p-1)p(p+1)}{12}$.

§3. Double Cosets and the Tits Building

We recall that $|T_3|$ is a 1-dimensional $\mathrm{SL}_3(\mathbb{Z})$–complex with 2 orbits of zero cells, and one orbit of 1-cells.* The respective isotropy subgroups are $P_1$, $P_2$ and $B$.

From the usual induction-restriction formula, we have

$$\mathbb{Z}[\mathrm{SL}_3(\mathbb{Z})/B]_{|\Gamma} \cong \bigoplus_{x_i \in I_0} \mathbb{Z}[\Gamma/\Gamma \cap x_i B x_i^{-1}]$$

where $I_0 = \Gamma\backslash\mathrm{SL}_3(\mathbb{Z})/B$ is the set of double cosets. Similarly we have

$$\mathbb{Z}[\mathrm{SL}_3(\mathbb{Z})/P_1]_{|\Gamma} \cong \bigoplus_{y_j \in I_1} \mathbb{Z}[\Gamma/\Gamma \cap y_j P_1 y_j^{-1}]$$

$$\mathbb{Z}[\mathrm{SL}_3(\mathbb{Z})/P_2]_{|\Gamma} \cong \bigoplus_{z_k \in I_2} \mathbb{Z}[\Gamma/\Gamma \cap z_k P_2 z_k^{-1}]$$

where $I_1 = \Gamma\backslash\mathrm{SL}_3(\mathbb{Z})/P_1$, $I_2 = \Gamma\backslash\mathrm{SL}_3(\mathbb{Z})/P_2$.

Note that as $\Gamma \lhd \mathrm{SL}_3(\mathbb{Z})$, we have $\Gamma \cap x_i B x_i^{-1} \cong \Gamma \cap B$, $\Gamma \cap y_j P_1 y_j^{-1} \cong \Gamma \cap P_1$, $\Gamma \cap z_k P_2 z_k^{-1} \cong \Gamma \cap P_2$, for all $i, j$ and $k$.

Next note that the cardinality of the indexing can be determined as follows*:

$$\#I_0 = [\mathrm{SL}_3(\mathbb{Z})/\Gamma : B/\Gamma \cap B] = [\mathrm{SL}_3(p) : B(p)]$$

$$\#I_1 = [\mathrm{SL}_3(\mathbb{Z})/\Gamma : P_1/\Gamma \cap P_1] = [\mathrm{SL}_3(p) : P_1(p)]$$

$$\#I_2 = [\mathrm{SL}_3(\mathbb{Z})/\Gamma : P_2/\Gamma \cap P_2] = [\mathrm{SL}_3(p) : P_2(p)]$$

* By the 'Invariant Factor Theorem', $\mathrm{SL}_3(\mathbb{Z})$ acts transitively on one and two dimensional subspaces in $\mathbb{Q}^3$, as well as on the flags.

* These indexing sets will also correspond to the number of $\Gamma$–conjugacy classes of parabolic subgroups.

where $B(p)$, $P_1(p)$, $P_2(p)$ denote the corresponding subgroups in $\mathrm{SL}_3(p)$.

The orders of these subgroups can easily be determined to be

$$|B(p)| = 4p^3 \quad \text{and} \quad |P_1(p)| = |P_2(p)| = 2p^3(p^2-1).$$

From this we readily deduce that

$$\#I_0 = \frac{(p^3-1)(p^2-1)}{4} \quad \text{and} \quad \#I_1 = \#I_2 = \frac{p^3-1}{2}.$$

We now use this to describe the cellular chains on $|T_3|$ as $\Gamma$-modules. Given that over $\mathrm{SL}_3(\mathbb{Z})$ we have

$$C_0(|T_3|) \cong \mathbb{Z}[\mathrm{SL}_3(\mathbb{Z})/P_1] \oplus \mathbb{Z}[\mathrm{SL}_3(\mathbb{Z})/P_2]$$

$$C_1(|T_3|) \cong \mathbb{Z}[\mathrm{SL}_3(\mathbb{Z})/B]$$

we obtain that as $\Gamma$-modules:

$$C_0(|T_3|) \cong \left(\bigoplus^{\frac{p^3-1}{2}} \mathbb{Z}[\Gamma/\Gamma \cap y_j P_1 y_j^{-1}]\right) \oplus \left(\bigoplus^{\frac{p^3-1}{2}} \mathbb{Z}[\Gamma/\Gamma \cap z_k P_2 z_k^{-1}]\right)$$

$$C_1(|T_3|) \cong \bigoplus^{\frac{(p^3-1)(p^2-1)}{4}} \mathbb{Z}[\Gamma/\Gamma \cap x_i B x_i^{-1}].$$

We note here that dividing out by the $\Gamma$-action we can describe $|T_3|/\Gamma$ in terms of the mod $p$ building associated to $\mathrm{SL}_3(p)$. This a $(p+1)$-valent graph which can be described explicitly.

## §4. Cohomology Calculations

We now use the second spectral sequence described in §0, applied to the 1-dimensional $\Gamma$-complex $|T_3|$; we obtain a long exact sequence

$$\cdots \to H^i_\Gamma(|T_3|) \to H^i(\Gamma, C^0(|T_3|)) \to H^i(\Gamma, C^0(|T_3|)) \to H^{i+1}_\Gamma(|T_3|) \to \cdots$$

Combining the results in §2, §3, we have

$$H^*(\Gamma, C^0(|T_3|)) \cong \left(\bigoplus^{\frac{p^3-1}{2}} H^*(P_1 \cap \Gamma(p), \mathbb{Z})\right) \oplus \left(\bigoplus^{\frac{p^3-1}{2}} H^*(P_2 \cap \Gamma(p), \mathbb{Z})\right)$$

$$H^*(\Gamma, C^1(|T_3|)) \cong \bigoplus^{\frac{(p^3-1)(p^2-1)}{4}} H^*(B \cap \Gamma(p), \mathbb{Z}).$$

We are now ready to substitute explicit values:

$$0 \to \mathbb{Z} \to \mathbb{Z}^{p^3-1} \to \mathbb{Z}^{\frac{(p^3-1)(p^2-1)}{4}} \to H^1_\Gamma(|T_3|) \to \mathbb{Z}^{\left[\frac{12+(p-1)p(p+1)}{6}\right]\left[\frac{p^3-1}{2}\right]} \to$$

$$\text{(4.1)}\quad \mathbb{Z}^{\frac{(p^3-1)(p^2-1)}{2}} \to H^2_\Gamma(|T_3|) \to \left[\bigoplus^{\frac{p^3-1}{2}} H^1(\Gamma_2(p), M\oplus M^*)\oplus\mathbb{Z}\oplus\mathbb{Z}\right] \to$$

$$\left[\bigoplus^{\frac{(p^3-1)(p^2-1)}{4}} \mathbb{Z}\oplus\mathbb{Z}\oplus\mathbb{Z}/p\right] \to H^3_\Gamma(|T_3|) \to \bigoplus^{p^3-1}\mathbb{Z}^{\left[\frac{12+(p-1)p(p+1)}{12}\right]} \to \mathbb{Z}^{\left[\frac{(p^3-1)(p^2-1)}{4}\right]} \to$$

$$H^4_\Gamma(|T_3|) \to 0.$$

We now derive some consequences from this sequence. Recall first that $\chi(B\Gamma(p)) = 0$ and hence from (1.2) we deduce $\chi(BQ(p)) = \chi(|T_3| \times_\Gamma E\Gamma) = 0$. Hence if $\mathbb{F}$ is a field of characteristic prime to $p$, we infer from (4.1) that

$$\dim_{\mathbb{F}} H^1(\Gamma_2(p), M_{\mathbb{F}} \oplus M^*_{\mathbb{F}}) = \frac{(p-1)p(p+1)}{3}$$

where $M_{\mathbb{F}}$ is the 2-dimensional natural representation of $\Gamma_2(p)$. In addition one can easily verify that if $(p,q)=1$, then $H^1(\Gamma_2(p), M\oplus M^*)$ is $q$-torsion-free.

At the prime $p$, we obtain a long exact sequence associated to mod $p$ reduction, namely

$$0 \to (\mathbb{F}_p)^4 \to H^1(\Gamma_2(p), M\oplus M^*) \xrightarrow{\cdot p} H^1(\Gamma_2(p), M\oplus M^*) \to [\mathbb{F}_p]^{4\left(1+\frac{(p-1)p(p+1)}{12}\right)} \to 0.$$

Observing that $[(M\oplus M^*)\otimes \mathbb{Z}/p^2]^{\Gamma_2(p)} \cong (\mathbb{Z}/p)^4$, we deduce

**Proposition 4.2**.

$$H^1(\Gamma_2(p), M\oplus M^*) \cong [\mathbb{Z}]^{\frac{(p-1)p(p+1)}{3}} \oplus (\mathbb{Z}/p^2)^4. \qquad \square$$

Let $\beta_3$ denote the third betti number associated to the classifying space of a group. Looking at the final part of the sequence (4.1), we find that

$$\beta_3(Q(p)) \geq (p^3-1)\left[1+\frac{(p-1)p(p+1)}{12}\right] - \frac{(p^3-1)(p^2-1)}{4} + 1$$

which, combined with (1.2) and the fact that $H^1(\Gamma,\mathbb{Z}) = 0$ yields

**Theorem 4.3**.

$$\beta_3(\Gamma(p)) \geq \beta_3(Q(p)) \geq \tfrac{1}{12}(p^3-1)(p^3-3p^2-p+15)+1. \qquad \square$$

## §5 The Symplectic Group: Parabolic Subgroups

Unlike the case of the congruence subgroups in $\mathrm{SL}_3(\mathbb{Z})$, it appears that not much is known about the cohomology of the congruence subgroups of the symplectic group $\mathrm{Sp}_4(\mathbb{Z})$. We will use the methods previously described to analyze the cohomology of these groups. The key fact here is that $\mathrm{Sp}_4(\mathbb{Z})$ is also a rank two group, hence its Tits Building will once again be a graph with an edge transitive action of the group. As before, the cohomology can be approached using the parabolic subgroups.

We proceed as before, now let $|W|$ denote the geometric realization of the Tits Building for $\mathrm{Sp}_4(\mathbb{Z})$; we denote $L(p) = \pi_1(|W| \times_{\Gamma(p)} E\Gamma(p))$, where now $\Gamma(p) \subset \mathrm{Sp}_4(\mathbb{Z})$ is a level $p$ congruence subgroup, $p$ an odd prime. In this instance, $\Gamma(p)$ has cohomological dimension equal to 4, and hence $L(p)$ has cohomological dimension equal to 5; as before it will be a Poincaré Duality group.

We need to describe the parabolic subgroups in some detail. If

$$J = \begin{pmatrix} 0 & 0 & 0 & -1 \\ 0 & 0 & -1 & 0 \\ 0 & 1 & 0 & 0 \\ 1 & 0 & 0 & 0 \end{pmatrix}$$

then $\mathrm{Sp}_4(\mathbb{Z})$ can be described as the matrices $A \in \mathrm{SL}_4(\mathbb{Z})$ such that $A^T J A = J$. It is convenient to write

$$A = \begin{pmatrix} A_1 & A_2 \\ A_3 & A_4 \end{pmatrix}$$

where $A_1, A_2, A_3, A_4$ are $2 \times 2$ matrices. Then, if we write

$$J = \begin{pmatrix} 0 & -Q \\ Q & 0 \end{pmatrix},$$

the conditions on $A$ can be rewritten as

$$A_3^T Q A_1 - A_1^T Q A_3 = 0, \;\; A_3 Q A_2 - A_1^T Q A_4 = -Q$$

$$A_4^T Q A_1 - A_2^T Q A_3 = Q, \;\; A_4^T Q A_2 - A_2^T Q A_4 = 0.$$

The parabolic subgroups in $\mathrm{Sp}_4(\mathbb{Z})$ can be obtained by intersecting with appropriate parabolics in $\mathrm{SL}_4(\mathbb{Z})$. To begin we consider the parabolic subgroup $G_1$ of symplectic matrices with $A_3 = 0$. Manipulating the defining equations one can easily see that if

$$S = \begin{pmatrix} I & 0 \\ 0 & Q \end{pmatrix},$$

then the conjugate group $SG_0S^{-1}$ can be described as the set of matrices

$$\begin{pmatrix} D & B \\ 0 & (D^T)^{-1} \end{pmatrix}$$

where $D \in \mathrm{GL}_2(\mathbb{Z})$ and $BD^T = DB^T$.

Next we consider the symplectic matrices where the first column has zero entries except in the top corner. To simplify things we will assume that this non–zero entry is 1. After some manipulation we obtain a description of this group (denoted $G_2$) as the matrices of the form

$$\begin{pmatrix} 1 & (u_1\ u_2)\begin{pmatrix} 0 & 1 \\ -1 & 0 \end{pmatrix} C & & w \\ 0 & c_{11} & c_{12} & u_1 \\ 0 & c_{21} & c_{22} & u_2 \\ 0 & 0 & 0 & 1 \end{pmatrix}$$

where

$$C = \begin{pmatrix} c_{11} & c_{12} \\ c_{21} & c_{22} \end{pmatrix} \in \mathrm{SL}_2(\mathbb{Z}),\ u_1, u_2, w \in \mathbb{Z}.$$

Letting $C \in \mathrm{GL}_2(\mathbb{Z})$ we obtain an index 2 extension group which is evidently a maximal parabolic, denoted $G_2$.

The intersection $G_0 = G_1 \cap G_2$ can then be described as the determinant one matrices of the form

$$\begin{pmatrix} \pm 1 & x & y & z \\ 0 & \pm 1 & w & y - xw \\ 0 & 0 & \pm 1 & -x \\ 0 & 0 & 0 & \pm 1 \end{pmatrix}$$

where $x, y, z, w \in \mathbb{Z}$. As expected, this coincides with the intersection $U_T \cap \mathrm{Sp}_4(\mathbb{Z})$, where $U_T$ is the subgroup of all $4 \times 4$ matrices which are upper triangular. We will not dwell on the parabolics themselves, rather we are interested in their intersections with the congruence subgroup $\Gamma(p)$.

**5.1 The group $G_1 \cap \Gamma(p)$:**

Using the fact that $\Gamma(p) \lhd \mathrm{Sp}_4(\mathbb{Z})$, this group can be described as the matrices of the form

$$\begin{pmatrix} D & B \\ 0 & (D^T)^{-1} \end{pmatrix}$$

where $D \in \Gamma_2(p) \subset \mathrm{SL}_2(\mathbb{Z})$, $B = pE$ and $ED^T = DE^T$. Let $V$ denote the subgroup consisting of all matrices in $G_1 \cap \Gamma(p)$ of the form

$$\begin{pmatrix} I & B \\ 0 & I \end{pmatrix};$$

which one can easily check is a normal subgroup isomorphic to $\mathbb{Z}^3$. Using this it follows that

$$G_0 \cap \Gamma(p) \cong V \times_T \Gamma_2(p)$$

where $X \in \Gamma_2(p)$ acts via

$$\begin{pmatrix} I & B \\ 0 & I \end{pmatrix} \mapsto \begin{pmatrix} I & XBX^T \\ 0 & I \end{pmatrix}.$$

In fact this action can be expressed more generally as the restriction of a homomorphism $T : \mathrm{SL}_2(\mathbb{Z}) \to \mathrm{SL}_3(\mathbb{Z})$ given by

$$\begin{pmatrix} a & b \\ c & d \end{pmatrix} \mapsto \begin{pmatrix} a^2 & 2ab & b^2 \\ ac & ad+bc & bd \\ c^2 & 2cd & d^2 \end{pmatrix}.$$

### 5.2 The group $G_2 \cap \Gamma(p)$:

In the description of $G_2$ given above, we simply require that $C \in \Gamma_2(p) \subset \mathrm{SL}_2(\mathbb{Z})$, $u_1, u_2, z \in p\mathbb{Z}$ and that the other 2 diagonal entries be equal to one. Let $U$ denote the normal subgroup consisting of matrices in $G_2 \cap \Gamma(p)$ such that $C = I$. Then we can write the group as a semidirect product

$$U \times_T \Gamma_2(p)$$

where an element $C \in \Gamma_2(p)$ acts via

$$\begin{pmatrix} 1 & (u_1\ u_2)\begin{pmatrix} 0 & 1 \\ -1 & 0 \end{pmatrix} & z \\ 0 & I & \begin{pmatrix} u_1 \\ u_2 \end{pmatrix} \\ 0 & 0 & 1 \end{pmatrix} \mapsto \begin{pmatrix} 1 & (u_1\ u_2)\begin{pmatrix} 0 & 1 \\ -1 & 0 \end{pmatrix} C^{-1} & z \\ 0 & I & C\begin{pmatrix} u_1 \\ u_2 \end{pmatrix} \\ 0 & 0 & 1 \end{pmatrix}$$

### 5.3 The group $G_0 \cap \Gamma(p)$:

We conjugate to get this into suitable form:

$$\begin{pmatrix} I & 0 \\ 0 & Q \end{pmatrix} \begin{pmatrix} 1 & x & y & z \\ 0 & 1 & w & y-xw \\ 0 & 0 & 1 & -x \\ 0 & 0 & 0 & 1 \end{pmatrix} \begin{pmatrix} I & 0 \\ 0 & Q \end{pmatrix} = \begin{pmatrix} 1 & x & z & y \\ 0 & 1 & y-xw & w \\ 0 & 0 & 1 & 0 \\ 0 & 0 & -x & 1 \end{pmatrix}$$

from which we deduce that

$$G_0 \cap \Gamma(p) \cong V \times_T \mathbb{Z}$$

where the group $\mathbb{Z}$ is generated by

$$t = \begin{pmatrix} 1 & p \\ 0 & 1 \end{pmatrix},\ T(t) = \begin{pmatrix} 1 & 2p & p^2 \\ 0 & 1 & p \\ 0 & 0 & 1 \end{pmatrix}.$$

## §6 The Symplectic Group: Cohomology Calculations

We are now in a position to calculate the cohomology of the groups $G_i \cap \Gamma(p)$. Using the spectral sequence associated to the extension in 5.1, we have

**Theorem 6.1**

$$H^*(G_1 \cap \Gamma(p), \mathbb{Z}) \cong \begin{cases} \mathbb{Z} & * = 0 \\ (\mathbb{Z})^{\frac{1+(p-1)p(p+1)}{12}} & * = 1 \\ H^1(\Gamma_2(p), H^1(V, \mathbb{Z})) & * = 2 \\ \mathbb{Z} \oplus H^1(\Gamma_2(p), H^2(V, \mathbb{Z})) & * = 3 \\ (\mathbb{Z})^{1+\frac{(p-1)p(p+1)}{12}} & * = 4. \\ 0 & \text{otherwise.} \end{cases}$$

□

To compute $H^*(G_2 \cap \Gamma(p), \mathbb{Z})$, we first observe that

$$H^*(U, \mathbb{Z}) \cong \begin{cases} \mathbb{Z} & * = 0 \\ \mathbb{Z} \oplus \mathbb{Z} & * = 1 \\ \mathbb{Z} \oplus \mathbb{Z} \oplus \mathbb{Z}/p & * = 2 \\ \mathbb{Z} & * = 3 \end{cases}.$$

As a $\Gamma_2(p)$–module, the cohomology of $H^*(U, \mathbb{Z})$ is trivial in dimensions 0 and 3; the dual of the 2-dimensional natural representation $M^*$ in $H^1$; and $M^* \oplus \mathbb{Z}/p$, where $\mathbb{Z}/p$ has a trivial action in $H^2$. From this information we can obtain

**Theorem 6.2**

$$H^*(G_2 \cap \Gamma(p), \mathbb{Z}) \cong \begin{cases} \mathbb{Z} & * = 0 \\ H^1(\Gamma_2(p), \mathbb{Z}) & * = 1 \\ H^1(\Gamma_2(p), M^*) \oplus \mathbb{Z}/p & * = 2 \\ \mathbb{Z} \oplus H^1(\Gamma_2(p), M^*) \oplus H^1(\Gamma_2(p), \mathbb{Z}/p) & * = 3 \\ (\mathbb{Z})^{1+\frac{(p-1)p(p+1)}{12}} & * = 4 \\ 0 & \text{otherwise.} \end{cases}$$

□

To compute $H^*(G_0 \cap \Gamma(p), \mathbb{Z})$ we use the semidirect product description provided in 5.3 to obtain

$$H^*(G_0 \cap \Gamma(p), \mathbb{Z}) \cong H^0(\mathbb{Z}, H^*(V, \mathbb{Z})) \oplus H^1(\mathbb{Z}, H^{*-1}(V, \mathbb{Z})).$$

These terms can be explicitly computed, yielding

**Theorem 6.3**

$$H^*(G_0 \cap \Gamma(p), \mathbb{Z}) \cong \begin{cases} \mathbb{Z} & * = 0 \\ \mathbb{Z} \oplus \mathbb{Z} & * = 1 \\ \mathbb{Z} \oplus \mathbb{Z} \oplus \mathbb{Z}/p \oplus \mathbb{Z}/2p & * = 2 \\ \mathbb{Z} \oplus \mathbb{Z} \oplus \mathbb{Z}/p \oplus \mathbb{Z}/2p & * = 3 \\ \mathbb{Z} & * = 4. \\ 0 & \text{otherwise.} \end{cases}$$

□

Just as before, the cohomology of these subgroups can be used to estimate the cohomology of $\Gamma(p)$. We will concentrate here on finding a lower bound for $\beta_4(\Gamma(p))$. For this we need to compute the indices $j_i = [\mathrm{Sp}_4(\mathbb{F}_p) : G_i(p)]$ where $G_i(p)$ denote the corresponding groups reduced mod $p$. The following can readily be established:

$$|\mathrm{Sp}_4(\mathbb{F}_p)| = p^4(p^4-1)(p^2-1), \ |G_0(p)| = 8p^4, \ |G_1(p)| = 2p^4(p^2-1) = |G_2(p)|$$

From this we deduce

$$j_0 = (p^4-1)(p^2-1)/8, \ j_1 = j_2 = (p^4-1)/2.$$

As in §4 we have a long exact sequence which can be used to estimate the cohomology of $\Gamma(p)$. In this case when we look at the end of the sequence we obtain a lower bound for $\beta_4(L(p))$ and hence we have

**Theorem 6.4**

$$\beta_4(\Gamma(p)) \geq \beta_4(L(p)) \geq \frac{1}{24}(p^4-1)(2p^3-3p^2-2p+27)+1.$$

□

Note that we have a complete *integral* calculation for the cohomology of the constituent subgroups. Hence a calculation over the field $\mathbb{F}_p$ is also possible. Also it is interesting to note the presence of 2–torsion in the cohomology of $G_0 \cap \Gamma(p)$.

§7 Final Remarks

Let $\mathcal{A} = P_1 *_B P_2$ denote the usual amalgamated product of parabolics in $\mathrm{SL}_3(\mathbb{Z})$. This group acts in the usual way on an infinite tree $K$ with a single edge as a quotient. Now there is a natural surjection $\mathcal{A} \to \mathrm{SL}_3(\mathbb{Z})$ induced by the inclusions of the parabolic subgroups. If the kernel of this map is denoted by $F$, then it is the fundamental group of an infinite graph on which $\mathrm{SL}_3(\mathbb{Z})$ acts with quotient a single edge. In fact this graph can be identified with the Tits building, and from our previous constructions it is easy to see that the group $Q(p)$ is a finite index subgroup fitting into an extension

$$1 \to Q(p) \to \mathcal{A} \to \mathrm{SL}_3(\mathbb{F}_p) \to 1.$$

In other words, the amalgam $\mathcal{A}$ is a 4-dimensional virtual duality group, mapping onto $\mathrm{SL}_3(\mathbb{Z})$. In fact this is a purely algebraic construction of the Borel–Serre compactification for $\mathrm{SL}_3(\mathbb{Z})$. More precisely, it is easy to verify that if $\overline{X}$ is their construction, then $\pi_1(\partial(\overline{X}/\Gamma(p))) \cong Q(p)$. Observe that in this case $\partial(\overline{X}/\Gamma(p))$ is a closed 4-dimensional manifold, in fact homotopy equivalent to the Borel construction $|T_3| \times_{\Gamma(p)} E\Gamma(p) \cong BQ(p)$. This situation evidently generalizes to other rank 2 situations, such as $\mathrm{Sp}_4(\mathbb{Z})$ (discussed in §6) and $\mathrm{G}_2(\mathbb{Z})$. It is hoped that the purely cohomological approach which we have presented here will help yield some insight for understanding the cohomology of congruence subgroups, but it remains a difficult undertaking.

Detailed calculations for the cohomology of congruence subgroups with rational coefficients can also be found in the paper by Lee and Schwermer [LS], we refer the interested reader to this reference for further details. In their language, the image of the cohomology of $\Gamma(p)$ in the cohomology of $Q(p)$ corresponds to the 'cohomology at infinity' and the kernel can be identified with part of the 'cuspidal cohomology'.